\documentclass[10pt,reqno]{amsart}
     \makeatletter
     \def\section{\@startsection{section}{1}%
     \z@{.7\linespacing\@plus\linespacing}{.5\linespacing}%
     {\bfseries
     \centering
     }}
     \def\@secnumfont{\bfseries}
     \makeatother
\newtheorem{theorem}{Theorem}[section]
\newtheorem{lemma}[theorem]{Lemma}

\theoremstyle{definition}
\newtheorem{definition}[theorem]{Definition}

\theoremstyle{remark}

\numberwithin{equation}{section} 

\usepackage{graphicx}
\usepackage{amssymb}
\usepackage{graphicx}

\begin{document}
\title[Percolation in a hierarchical random graph]{Percolation in a hierarchical random graph}

\author[D.~A.~Dawson]{D.~A.~Dawson*}
\address{D.A. Dawson: School of Mathematics and Statistics, Ottawa, Canada K1S 5B6}
\email{ddawson@math.carleton.ca}
\thanks{* This research is  supported by NSERC}

\author[L.G. Gorostiza]{L.~G.~Gorostiza**}
\thanks{** Partially supported by CONACyT grant 45684-F}
\address{L.G. Gorostiza: CINVESTAV, Mexico City, Mexico}
\email{lgorosti@math.cinvestav.mx}



\begin{abstract}
 We study asymptotic
percolation as $N\to \infty$ in an infinite random graph
${\mathcal G}_N$ embedded in the hierarchical group of order $N$,
with connection probabilities depending on an ultrametric distance
between vertices. ${\mathcal G}_N$ is structured as a cascade of
finite random subgraphs of (approximate) Erd$''\kern -.20cm{\rm
o}$s-R\'enyi type. However, the results are different from those
of classical random graphs, e.g.,  the average length of paths in
the giant component of an ultrametric ball is much longer than in
the classical case. We give a criterion for percolation, and show
that percolation takes place along giant components of giant
components at the previous level in the cascade of  subgraphs for
all consecutive hierarchical distances. The proof involves a
hierarchy of ``doubly stochastic''  random graphs with vertices
having an internal structure and random connection probabilities.

\end{abstract}

\maketitle

\section{Introduction}
\setcounter{section}{1} \setcounter{equation}{0} \label{sec:1}
Random graphs have been used to analyze percolation in some
infinite systems (e.g., \cite{BLP,BS} and references therein). On
the other hand, hierarchical structures have been used in
applications in physics,  biology (in particular, genetics), etc.,
where an underlying ultrametric distance plays a basic role (e.g.,
\cite{CDG, DGW1, DGW2, SF} and references therein). Hence it is
natural to consider infinite random graphs embedded in
hierarchical structures, with the probability of connection
between two vertices depending on an ultrametric distance between
them. In this paper we study percolation in an infinite random
graph, ${\mathcal G}_N$, embedded in an ultrametric group,
$\Omega_N$, called hierarchical group of order $N$, as
$N\to\infty$. The structure of ${\mathcal G}_N$ is a cascade of
(approximate) Erd$''\kern -.20cm{\rm o}$s-R\'enyi random subgraphs
at consecutive hierarchical distances, and this allows using
results from the classical theory of random graphs \cite{AS,B,D,
ER, JLR}. However, we use them only as a technical tool. As we
shall see, for ${\mathcal G}_N$ there   are  different results
from the classical ones. An essential feature of ${\mathcal G}_N$
is that, due to the tree-type form of $\Omega_N$, any path
connecting two vertices must contain an edge (or 1-step
connection) of size equal to the hierarchical distance between
them. In order to prove that percolation occurs along giant
components in the cascade at consecutive distances, we introduce a
family of random graphs whose vertices have an internal structure
involving giant components at the previous hierarchical levels. In
these ``doubly stochastic'' random graphs  the connection
probabilities between vertices are  random variables which are
highly dependent, but asymptotically deterministic as
$N\rightarrow \infty$. There are papers that consider random
graphs with connection probabilities depending on a distance
between vertices or some function of pairs of vertices (e.g.,
\cite{BB, Bi, BJR, CG}),
 and/or involving hierarchical structures or ultrametrics (e.g., \cite{BR, Bi, HJ, H, K1, K2,  RSM, VCY}), but
  we have not seen the setup we consider here.

\begin{definition}
\label{D1.1} {\rm The {\it hierarchical group of order $N$} (integer $\geq 2$) is defined as
\begin{eqnarray*}
 \Omega_N\!&=&\!\{{\mathbf x}=(x_1, x_2, \ldots):x_i \in
\{0,1,\ldots, N-1\},\; i=1,2,\ldots, ;\\&&
\qquad\qquad\qquad\qquad\qquad\quad x_i \neq 0\;\hbox{\rm only for
finitely many}\;i\},
\end{eqnarray*}
with addition componentwise mod $N$, i.e., $\Omega_N$ is the direct sum of a countable
number of copies of the cyclic group of order $N$. The {\it hierarchical distance} on $\Omega_N$ is given by
\begin{eqnarray*}
d({\mathbf x},{\mathbf y})=\left\{
\begin{array}{lcc}
0 & if & {\mathbf x}={\mathbf y},\\
\max \{i:x_i\neq y_i \} & if & {\mathbf x}\neq {\mathbf y}.
\end{array}\right.
\end{eqnarray*}}
\end{definition}

Note that $d(\cdot, \cdot)$ is translation-invariant, and it is  an ultrametric, i.e., it satisfies the strong
(or non-Archimedean) triangle inequality: for any ${\mathbf x},{\mathbf y},{\mathbf z}$,
\[
d({\mathbf x},{\mathbf y}) \leq \max \{d({\mathbf x},{\mathbf z}), d({\mathbf z},{\mathbf y})\},
\]
and for each integer $k>1$, a {\it $k$-ball} in $\Omega_N$ (i.e., a set of points at distance at most $k$ from each other) contains $N^k$ points, and it is the union of $N$  $(k-1)$-balls which are at distance $k$ from each other
(a $0$-ball is a single point).

The group $\Omega_N$ has been used as a state space for stochastic models in several applications (e.g., \cite{DGW1} and references therein).

We are interested in studying percolation properties on
$\Omega_N$, and to this end we consider a random graph ${\mathcal
G}_N$ whose vertex set is $\Omega_N$, with connection
probabilities depending on the hierarchical distance between
points. We parameterize the connection probabilities so as to
characterize the critical regime for percolation.

The random graph ${\mathcal G}_N$ is defined as follows.

\begin{definition}
\label{D1.2}{\rm We define an infinite random graph ${\mathcal
G}_N$ with the points of $\Omega_N$ as vertices,
 and for each $k\geq 1$ the probability of connection between
${\mathbf x}$ and ${\mathbf y}$ such that $d({\mathbf x},{\mathbf y})=k$ is given by
\begin{equation}
\label{eq:1.4}
\frac{c_k}{N^{2k-1}},
\end{equation}
where $c_k$ is a positive constant (independent of $N$), all connections being independent.

In this paper we do not investigate percolation  for fixed $N$, we
obtain a result  for large $N$, i.e., asymptotic percolation as
$N\to \infty$, by means of the theory of random graphs.   This
method is natural for the framework of the hierarchy of  random
graphs which constitutes ${\mathcal G}_N$.

Recall that in the classical theory of Erd$''\kern -.20cm{\rm
o}$s-R\'enyi ({\it ER}) random graphs ${\mathcal G}(N,p)$ with $N$
vertices and connection probability $p$, the right
parameterization for existence of giant components as
$N\rightarrow \infty$ is $p=c/N$, $c>1$. Analogously, the choice
of connection probabilities (\ref{eq:1.4}) will allow us to
characterize percolation in ${\mathcal G}_N$ as $N\rightarrow
\infty$ in terms of the sequence $(c_k)$ with all $c_k >1$. One of
our main tools is the {\it ER} theory, which we will apply
recursively along consecutive hierarchical distances. Since a
$k$-ball has $N^k$ points, the {\it ER} theory would require
connection probabilities $c_k /N^k, c_k>1$. The fact that
percolation in ${\mathcal G}_N$ can occur with connection
probabilities of the form (\ref{eq:1.4}) is due to the internal
hierarchical structure of $k$-balls. We shall see in Section 3.2
that the average length of shortest paths between randomly chosen
points in the giant component of a $k$-ball is of order $(\log
N)^k$, which should be compared with $k\log N$ in the classical
case.

To begin, we observe why the  theory of {\it ER} random graphs
can be used to investigate  ${\mathcal G}_N$. By (\ref{eq:1.4}),
two different $(k-1)$-balls in a given $k$-ball  are connected
within the $k$-ball  with probability
\begin{equation}
\label{eq:1.5}
p^{N,k}=1-\biggl(1-\frac{c_k}{N^{2k-1}}\biggr)^{N^{2(k-1)}}, \quad k>1,\quad p^{N,1}=\frac{c_1}{N}.
\end{equation}
Hence $p^{N,k} \sim c_k/N$ as $N\rightarrow \infty$. More precisely, using the elementary inequality
\begin{equation}
\label{eq:1.6}
0<\frac{m}{n}y -\biggl[1-\biggl(1-\frac{y}{n}\biggr)^m\biggr] <\biggl(\frac{m}{n}\biggr)^2\frac{y^2}{2},\quad 0<y<n, \quad m\geq 2,
\end{equation}
with $y=c_k$, $n=N^{2k-1}$, $m=N^{2(k-1)}$, we have from
(\ref{eq:1.5})
\begin{equation}
\label{eq:1.7}
p^{N,k} =\frac{c_k}{N}+o \biggl(\frac{1}{N}\biggr)\quad {\rm as} \quad
N\rightarrow \infty,\quad k>1, \quad p^{N,1}=\frac{c_1}{N}.
\end{equation}
Two $(k-1)$-balls in a $k$-ball can also be connected through points outside the $k$-ball. However, as
 we shall see, this does not add more than $o(1/N)$ to the probability of connection between the $(k-1)$-balls.
  On the other hand, we are interested in what happens inside the $k$-ball. Therefore we will disregard connections
  through points outside the $k$-ball. By (\ref{eq:1.7}), for any $k>1$,
a $k$-ball with its $N\,\,(k-1)$-balls as vertices is a random
graph ${\mathcal G}(N,p^{N,k})$ which approximates
 the {\it ER} random graph ${\mathcal G}(N,c_k/N)$ as
$N\to\infty$. Note that (\ref{eq:1.4}) is equal to $c_k/N$ with a norming which is the product of the
sizes of two $(k-1)$-balls. (See Remark 2.4).

The tree representation of $\Omega_N$ ({See Figure 1}) shows
${\mathcal G}_N$ as a cascade of random subgraphs ${\mathcal
G}(N,p^{N,k}), \;k=1,2,\dots$, where the form (1.4) of the
connection probabilities is inherited throughout the hierarchy.
Note that the hierarchical distances $k$ are the same for all $N$.


\begin{figure}
\begin{center}
\includegraphics[scale=0.75]{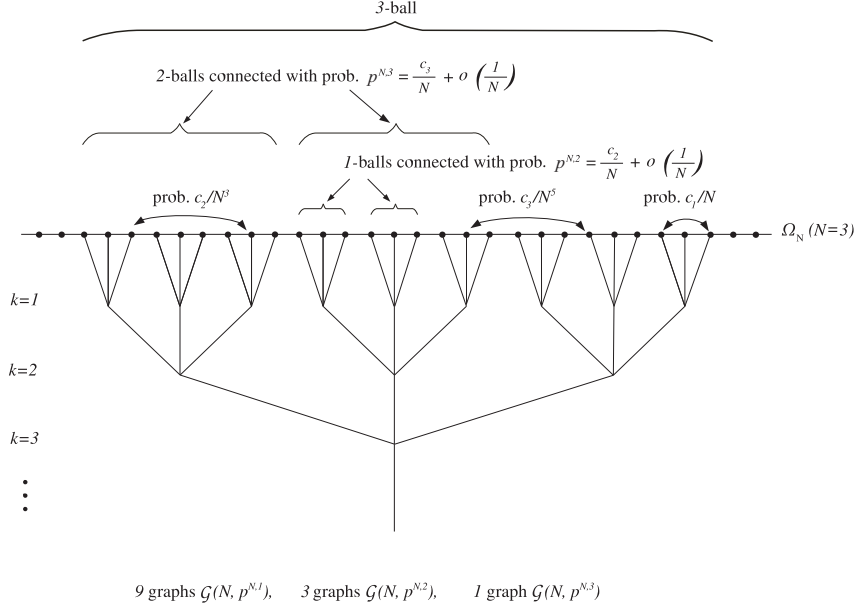}
\caption{Tree representation of $\Omega_N$}
\end{center}
\end{figure}

We consider asymptotic percolation in ${\mathcal G}_N$, meaning
that, as $N\rightarrow \infty$, there exists a path going from
${\mathbf 0}\in \Omega _N$ to infinity, i.e., to points at
arbitrarily large hierarchical distances, with positive
probability.  This is made precise in the following definition.}

\end{definition}
\begin{definition}\label{D1.3}{\rm
We say that there is  {\it asymptotic percolation}  in ${\mathcal
G}_N$ if
\begin{eqnarray}
\label{eq:1.8} &&P_{{\rm perc}}:=\inf_k\,\, \liminf_{N\rightarrow
\infty} P[{\mathbf 0}\in \Omega_N\;\; \hbox{\rm is connected by a
path to a point at distance}\;\; k]\nonumber\\&&>0.
\end{eqnarray}}
\end{definition}

In Theorem \ref{T2.1} we give a criterion  for asymptotic
percolation in ${\mathcal G}_N$, and we show that percolation
occurs along giant components of giant components at the previous
level in the cascade of  subgraphs ${\mathcal G}(N,p^{N,k})$ for
all hierarchical distances $k$. This involves replacing the
connection probabilities $p^{N,k}$, given by (\ref{eq:1.5}), by
random connection
 probabilities between giant components, which generates a hierarchy of doubly stochastic random graphs. Although
  we use  the basic results on the sizes of components in $ER$ random graphs, most of the work in the proof of the
  theorem consists in formulating the problem in such a way that those results can be employed to obtain upper and
  lower bounds for the probability of percolation. Another part of the proof consists in showing that percolation
   occurs along the cascade of giant components at consecutive hierarchical levels. We will  also discuss
 the possible form of percolation paths, and comment on approximate scale-free
degree sequences in the graphs ${\mathcal G}(N,c_k/N)$ for large
$k$.

In a final section we make some heuristic comments regarding doubly stochastic random graphs of the type used
for the proof of Theorem \ref{T2.1} and their application, in particular referring to average distance and central
limit theorem for the cascade of giant components. These random graphs give rise to questions that suggest further
research  which would be of independent interest.
\section{Asymptotic percolation and cascade percolation}
\setcounter{section}{2} \setcounter{equation}{0} We will  use the
following fundamental result on {\it ER} random graphs ${\mathcal
G}(N,c/N)$ (see, e.g., \cite{JLR}, Theorem 5.4; recall that
``a.a.s.'' for a graph property means that the probability that
the random graph
  possesses the property tends to $1$ as $N\rightarrow \infty$ \cite{JLR}):

\noindent (A) If $c>1$, then a.a.s. ${\mathcal G}(N,c/N)$ has  a
unique giant component, its size is  $\beta N$, where $\beta\in
(0,1)$ satisfies  $\beta=1-e^{-c\beta}$, and all the other
components have sizes  at most $\frac{16c}{(c-1)^2}\log N$.

 For each $k\geq 1$ and each $k$-ball, we
have a random graph ${\mathcal G}(N, p^{N,k})$
 whose vertices are the $N\;\;(k-1)$-balls contained in the $k$-ball, with connection probability
  $p^{N,k}$ given by (\ref{eq:1.5}). We assume that $c_k>1$ for all $k$, hence, because of (\ref{eq:1.7}) and (A),
   each ${\mathcal G}(N, p^{N,k})$ has a unique giant component a.a.s..
We may assume that $N$ is large enough so that
 the giant components have emerged in all $k$-balls for all $k$.
 Note that, by (\ref{eq:1.6}),  $p^{N,k}>q^{N,k}/N$, where $q^{N,k}=c_k (1-c_k/2N)>1$ for $N>c^2_k/2(c_k-1)$,
 therefore unique giant components emerge in ${\mathcal G}(N, p^{N,k})$
 for all $k$ as $N\rightarrow \infty$. Since percolation
 involves $c_k\to\infty$ as $k\to\infty$, our method is restricted
 to the limiting situation $N\to\infty$.

If $S$ is a non-empty set of vertices of ${\mathcal
G}(N,p^{N,k})$, we denote by $||S||$ the number of its elements
(hence $0<||S||\leq N$), and by $|S|$ the number of points of
$\Omega_N$ contained in $S$, thus $|S|=||S||N^{k-1}$.

It is reasonable to consider percolation with connections along
giant components in the  subgraphs ${\mathcal G}(N,p^{N,k})$ at
consecutive hierarchical distances $k$. Indeed, it may happen that
${\mathbf 0}$ is in the giant component of the $1$-ball it belongs
to, and that this $1$-ball is in the giant component of the
$2$-ball it belongs to, but it may also happen  that ${\mathbf 0}$
is not in the giant component of the  $1$-ball it belongs to, but
that it is connected to the giant component of the $2$-ball, which
may still be good for percolation. However, by (A), the second
possibility is negligible with respect to the first one as
$N\rightarrow \infty$. The same argument can be made for the
vertices of the subgraphs ${\mathcal G}(N,p^{N,k})$ and any two
consecutive values of $k$.
 In addition, if $Y^{(k)}$ denotes the  number of 1-step connections from a point to points at distance
  $k$ from it, we shall see that $P[Y^{(k)}>0]\sim c_k N^{-(k-1)}$ as $N\to\infty$ for each $k\geq 1$, so, the probability
  of connections at distance $k$ is negligible w.r.t. the probability of connections at distance $k-1$, and  all
  external connections from a $k$-ball will be at distance $k+1$
  from it for all $k$  as $N\to\infty$.
Therefore we will look at  asymptotic  percolation along giant components at
consecutive distances $k$. Moreover, we will consider  asymptotic percolation along the cascade of giant components
of giant components at the previous level for all hierarchical distances $k$.
 This precisely means that ${\mathbf 0}$ is in the giant component of the 1-ball it belongs
  to; this giant component is connected to other giant components of
$1$-balls in the $2$-ball in they belong to, forming a level $2$ giant component; and so
  on for every hierarchical distance $k$. For brevity, we call {\it cascade percolation} this special form of percolation.

In order to study cascade percolation, for each $k>1$ and large
$N$ we consider a random graph ${\mathcal G}_{N,k}$ whose vertices
are the  $(k-1)$-balls  that are the vertices of ${\mathcal G}(N,
p^{N,k})$, but now each vertex has an internal structure which
involves all  the components in the $j$-balls, $j=1,\ldots, k-2$,
in the $(k-1)$-ball. Clearly, the giant components are the ones
that will determine percolation. These random graphs are described
in the proof of Theorem 2.2, and their composition is reformulated
after the proof in terms of giant components alone.

For each $k\geq 1$, let $\beta_k \in (0,1)$ satisfy
\begin{equation}
\label{eq:2.2}
\beta_k =1-e^{-c_k \beta^2_{k-1}\beta_k}, \quad \beta_0=1,
\end{equation}
where $c_k\beta^2_{k-1}>1$. We will show that the probability of asymptotic percolation $P_{{\rm perc}}$
is given by $\prod \beta_k$. Hence we need a condition for
positivity of this product  in terms of the $c_k$, which are the
data of the graph ${\mathcal G}_N$. This is given in the following
lemma (proved in the Appendix).
\begin{lemma}
\label{L2.1}
Assume $c_k \nearrow \infty$ as $k\rightarrow \infty$, $c_1>2 \log 2$ and $c_2>8\log 2$. Then
\begin{equation}
\label{eq:2.3}
\beta_k >1/2 \quad \hbox{\it and}\quad c_k \beta^2_{k-1}>1\quad\hbox{\it  for all}\quad k,
\end{equation}

 \noindent and
\begin{equation}
\label{eq:2.4}
\prod^\infty_{k=1} \beta_k >0\quad\hbox{\it if and only if}\quad \sum^\infty_{k=1}e^{-c_k} <\infty.
\end{equation}
\end{lemma}

The main result in this paper is the following theorem.

\begin{theorem}
\label{T2.1} Assume $c_k \nearrow \infty$ as $k\rightarrow
\infty$, $c_1>2 \log 2$, $c_2>8\log 2$. Then there is asymptotic
percolation in ${\mathcal G}_N$ if and only if
$\sum^\infty_{k=1}e^{-c_k}<\infty$, asymptotic percolation occurs
in the form of cascade percolation, and the probability of
percolation  (\ref{eq:1.8}) is given by

\begin{equation}
\label{eq:2.5}
P_{{\rm perc}}=\prod^\infty_{k=1}\beta_k.
\end{equation}
\end{theorem}
\noindent
{\bf Proof.} Let us assume that percolation takes places only with connections along consecutive pairs of hierarchical distances $k$. In the last part of the proof we will show that this is so.

We prove first that  $\sum e^{-c_k}=\infty$ implies there is no
percolation in ${\mathcal G}_N$ along consecutive hierarchical
distances.
 By the left-hand inequality in (\ref{eq:1.6}), the connection probability $p^{N,k}$ given by (\ref{eq:1.5}) satisfies $$
p^{N,k}<\frac{c_k}{N}.
$$
Hence it suffices to show that there is no percolation along the
random graphs ${\mathcal G} (N,c_k/N)$ at distance $k$ for  every
$k$  in the weaker sense that ${\mathbf 0}$ is in the giant
component of the $1$-ball it belongs to,  this $1$-ball is in the
giant component of the
 $2$-ball it belongs to, and so on. The probability of percolation along the  graphs ${\mathcal G}(N, c_k/N)$
 is given by $\prod^\infty_{k=1}\gamma_k$  (connections at different distances are independent), where, by
 (A), $\gamma_k>0$ satisfies
\begin{equation}
\label{eq:2.6}
\gamma_k=1-e^{-c_k \gamma_k}
\end{equation}
for each $k\geq 1$. Now,  $\sum e^{-c_k}=\infty$ implies $\sum e^{-c_k \gamma_k}=\infty$,
since all $\gamma_k <1$, which by (\ref{eq:2.6}) implies $\sum (1-\gamma_k)=\infty$, and therefore $\prod \gamma_k=0$. Hence there is no  percolation. So, percolation along consecutive hierarchical distances implies $\sum e^{-c_k}<\infty$.

Note that if for some $k$  the vertex ($(k-1)$-ball) in ${\mathcal
G}(N,c_k/N)$ containing ${\mathbf 0}$ is not in the giant
component of the $k$-ball, then by (A) the corresponding
$\gamma_k$ in $\prod\gamma_k$ is $0$. More precisely, if the
$(k-1)$-ball containing ${\mathbf 0}$ is not in the giant
component in the $k$-ball, then, by (A), ${\mathbf 0}$ is
connected to at most $O (N^{k-1}\log N)$ points within the
$k$-ball. But the probability of connection between a set of
$N^{k-1}\log N$ points in the $k$-ball and points in the
$(k+1)$-ball outside the $k$-ball is, by (\ref{eq:1.4}),
\begin{eqnarray*}
\lefteqn{1-\biggl(1-\frac{c_{k+1}}{N^{2k+1}}\biggr)^{(N^{k-1}\log N)N^k(N-1)}}\\
&&\sim c_{k+1}\frac{\log N}{N}\rightarrow 0\quad {\rm as}\quad N\rightarrow \infty.
\end{eqnarray*}
Therefore there is no
percolation. This observation allows us to omit non-giant components in the proofs.

Now we assume $\sum e^{-c_k}<\infty$, hence $\prod \beta_k >0$ by Lemma \ref{L2.1}. We denote by $Q_{\rm perc}$ the probability of cascade percolation. We will show that
\begin{equation}
\label{eq:2.7'}
Q_{{\rm perc}}=\prod^\infty_{k=1}\beta_k.
\end{equation}
This obviously implies that there is asymptotic percolation in
${\mathcal G}_N$.

We prove first that
\begin{equation}
\label{eq:2.7N}
Q_{\rm perc}\geq \prod^\infty_{k=1}\beta_k.
\end{equation}
We proceed by steps.

For $k=2$, a given $2$-ball and large $N$, the graph ${\mathcal
G}_{N,2}$ has $N$ vertices, which are the $1$-balls in the
$2$-ball, and the internal structure of a $1$-ball is the family
of its components. Although only the giant component $G^{(1)}$
will be relevant for percolation, the other ones complete the size
of the $1$-ball, and this plays a role in the proof. By
(\ref{eq:1.4}), two giant components $G^{(1)}_i, G^{(1)}_j, i\neq
j$ (which are at distance  $2$ from each other) are connected with
conditional probability (given the sizes $|G^{(1)}_i|$,
$|G^{(1)}_j|$)
\begin{equation}
\label{eq:2.8N}
p^{N,2}_{ij}=1-\biggl(1-\frac{c_2}{N^3}\biggr)^{|G^{(1)}_i||G^{(1)}_j|}.
\end{equation}
We write $p^{N,2}_{ij}$ as
\begin{equation}
\label{eq:2.9N}
p^{N,2}_{ij}=\frac{c_2\beta^2_1}{N}+A^{N,2}_{ij}+B^{N,2}_{ij},
\end{equation}
where
\begin{eqnarray}
\label{eq:2.10N}
A^{N,2}_{ij}&=&1-\biggl(1-\frac{c_2}{N^3}\biggr)^{|G^{(1)}_i||G^{(1)}_j|}-\frac{|G^{(1)}_i||G^{(1)}_j|}{N^3}c_2,\\
\label{eq:2.11N}
B^{N,2}_{ij}&=&\biggl(\frac{|G^{(1)}_i ||G^{(1)}_j|}{N^3}-\frac{\beta^2_1}{N}\biggr)c_2.
\end{eqnarray}

By (\ref{eq:1.6}), $A^{N,2}_{ij}\leq 0$ and
$$
|A^{N,2}_{ij}|\leq \biggl(\frac{|G^{(1)}_j||G^{(1)}_j|}{N^3}\biggr)^2 \frac{c^2_2}{2}\leq \frac{c^2_2}{2N^2},
$$
since $|G^{(1)}_i|\leq N$. Hence
\begin{equation}
\label{eq:2.12N}
-\frac{c^2_2}{2N^2}\leq A^{N,2}_{ij}\leq 0.
\end{equation}

We take $0<\delta_2<1$ and delete from ${\mathcal G}_{N,2}$ the
vertices ($1$-balls) with giant components $G^{(1)}_i$ such that
\begin{equation}
\label{eq:2.13N}
|G^{(1)}_i|\leq \beta _1(1-\delta_2)N,
\end{equation}
 where $\beta_1$ satisfies (\ref{eq:2.2}) with $k=1$, except the one that contains ${\mathbf 0}$ if it happens to be in this case. For a pair $ij$ of giant components that do not satisfy (\ref{eq:2.13N}) we have, by (\ref{eq:2.11N}) and (\ref{eq:2.13N}),
\begin{equation}
\label{eq:2.14N}
B^{N,2}_{ij}>\frac{c_2\beta^2_1}{N}((1-\delta_2)^2-1)=\frac{c_2\beta^2_1}{N}(\delta^2_2-2\delta_2).
\end{equation}
Then, from (\ref{eq:2.9N}), (\ref{eq:2.12N}) and (\ref{eq:2.14N}), for the remaining giant components we have
\begin{equation}
\label{eq:2.15N}
p^{N,2}_{ij}>\frac{c_2\beta^2_1}{N}\biggl((1-\delta_2)^2-\frac{c_2}{2\beta^2_1N}\biggr).
\end{equation}

Since $c_2\beta^2_1>1$, by Lemma 2.1, for any $0<\varepsilon_2 <1$ such that $c_2\beta^2_1(1-\varepsilon_2)>1$ we can take $\delta_2$ small enough and $N_2$ large enough so that, from (\ref{eq:2.15N}),
\begin{equation}
\label{eq:2.16N}
p^{N,2}_{ij}>\frac{q^{(\varepsilon_2)}_2}{N},\quad \hbox{\rm where}\quad q^{(\varepsilon_2)}_2=c_2 \beta^2_1(1-\varepsilon_2)>1,
\end{equation}
for all $N\geq N_2$, for all pairs $ij$ of giant components which do not satisfy (\ref{eq:2.13N}).

Let $M^{N,2}$ denote the number of $1$-balls of ${\mathcal
G}_{N,2}$ whose giant components satisfy (\ref{eq:2.13N}).
$M^{N,2}$ is distributed  Bin$(N, P[|G^{(1)}|\leq \beta_1
(1-\delta_2)N])$ (the $|G^{(1)}_i|$ are i.i.d.). Hence, for any
$0<\eta <1$,
\begin{eqnarray*}
P\biggl[\frac{M^{N,2}}{N}>\eta \biggr] &\leq & \frac{1}{N\eta} EM^{N,2}=
\frac{1}{\eta}P[|G^{(1)}|\leq \beta_1(1-\delta_2)N]\\
&&\rightarrow 0\quad \hbox{\rm as}\quad N\rightarrow \infty,
\end{eqnarray*}
by  (\ref{eq:2.2}) for $k=1$ and (A) for the {\it ER} graph
${\mathcal G}(N, c_1/N)$ (see also \cite{AS}, p.167). So,
$M^{N,2}=o(N)$ as $N\rightarrow \infty$ ($o_p(N)$ in the notation
of \cite{JLR}, p.11). Therefore,  deleting from ${\mathcal
G}_{N,2}$ the $1$-balls whose giant components satisfy
(\ref{eq:2.13N}) does not impair the emergence of a giant
component in ${\mathcal G}_{N,2}$, which actually emerges due to
(\ref{eq:2.16N}). The giant component of ${\mathcal G}_{N,2}$ is a
{\it level 2} giant component, whose elements are $1$-balls with
their internal structures.

For $k=3$, a given $3$-ball and large $N$, the graph ${\mathcal
G}_{N,3}$ has $N$ vertices, which are the $2$-balls in the
$3$-ball, and the internal structure of a $2$-ball is the family
of its components, in particular the giant component of the
corresponding graph ${\mathcal G}_{N,2}$, as described above.
Having deleted from each $2$-ball an $o(N)$ number of $1$-balls as
above, we are left only with $1$-balls whose connection
probabilities between their giant components $p^{N,2}_{ij}$
satisfy (\ref{eq:2.16N}) (and the $1$-ball containing ${\mathbf
0}$). We now connect these $1$-balls with the smaller probability
$q^{(\varepsilon_2)}_2/N$ (recall that we are looking for a lower
bound for $Q_{{\rm perc}}$), and we consider the resulting level 2
giant components $G^{(2)}_i$.

By (\ref{eq:1.4}), two giant components $G^{(2)}_i, G^{(2)}_j, i\neq j$ (which are at distance $3$ from each other) are connected with probability
\begin{equation}
\label{eq:2.17'}
p^{N,3}_{ij}=1-\biggl(1-\frac{c_3}{N^5}\biggr)^{|G^{(2)}_i||G^{(2)}_j|},
\end{equation}
where $|G^{(2)}_i|=||G^{(2)}_i||N$ (each $1$-ball in $G^{(2)}_i$ contains $N$ points), and, by (A),
\begin{equation}
\label{eq:2.18'}
||G^{(2)}_i||\sim \beta^{(\varepsilon_2)}_2N\quad {\rm a.a.s.},
\end{equation}
where $\beta^{(\varepsilon_2)}_2>0$ satisfies
\begin{equation}
\label{eq:2.19'}
\beta^{(\varepsilon_2)}_2=1-e^{-q^{(\varepsilon_2)}_2\beta^{(\varepsilon_2)}_2}.\end{equation}

Now we are in a similar situation as in the previous step $(k=2)$,
with  $1$-balls playing the role of  points and
$q^{(\varepsilon_2)}_2$ playing the role of $c_1$. The fact that
$o(N)$ $1$-balls have been deleted from each $2$-ball has no
effect, and in any case it helps towards a lower bound. So,
proceeding as above with (\ref{eq:2.17'}), (\ref{eq:2.18'}),
(\ref{eq:2.19'}), we take $0<\delta_3<1$ and delete from
${\mathcal G}_{N,3}$ the $2$-balls with  giant components
$G^{(2)}_i$ such that
\begin{equation}
\label{eq:2.17N}
||G^{(2)}_i||\leq \beta^{(\varepsilon_2)}_2 (1-\delta_3)N,
\end{equation}
(except the one containing ${\mathbf 0}$). Choosing $\delta_3$ small enough and $N_3$ large enough (and larger than $N_2$),  we have that the connection probabilities $p^{N,3}_{ij}$ for the remaining giant components  satisfy
\begin{equation}
\label{eq:2.19N}
p^{N,3}_{ij}>\frac{q^{(\varepsilon_2, \varepsilon_3)}_3}{N},\quad {\rm where}\quad q^{(\varepsilon_2, \varepsilon_3)}_3=c_3(\beta^{(\varepsilon_2)}_2)^2 (1-\varepsilon_3)>1,
\end{equation}
with small enough $\varepsilon_3$, for all $N\geq N_3$. Note that since $q^{(\varepsilon_2)}_2\nearrow c_2\beta^2_1$ as $\varepsilon_2\searrow 0$, then, form (\ref{eq:2.19'}), $\beta^{(\varepsilon_2)}_2 \nearrow \beta_2$ as $\varepsilon_2 \searrow 0$, where $\beta_2$ satisfies (\ref{eq:2.2}) with $k=2$, and since $c_3\beta^2_2>1$ by Lemma 2.1, then $\varepsilon_2$ and $\varepsilon_3$ can be taken small enough so that $q^{(\varepsilon_2, \varepsilon_3)}_3>1$.

Since $\beta^{(\varepsilon_2)}_2<\beta_2<\gamma_2$ (see
(\ref{eq:2.6})), we  have as above that the number of deleted
$2$-balls from ${\mathcal G}_{N,3}$ is $o(N)$ as $N\rightarrow
\infty$, and we connect the remaining ones with the smaller
probability $q_3^{(\varepsilon_1, \varepsilon_3)}/N$.

Iterating this scheme, for each $k$, each given $k$-ball and large
$N$, we have a graph ${\mathcal G}_{N,k}$ whose $N$ vertices are
the $(k-1)$-balls in the $k$-ball with their internal structures,
which involve the level ($k-1$) giant components in the
corresponding graphs ${\mathcal G}_{N,k-1}$. A number $o(N)$ of
vertices of ${\mathcal G}_{N,k}$ are deleted, as above,  the
remaining ones are connected with the smaller probability
$q^{(\varepsilon_2, \ldots , \varepsilon_k)}_k/N$, where
\begin{equation}
\label{eq:2.20N}
q^{(\varepsilon_2,\ldots, \varepsilon_k)}_k =c_k (\beta^{(\varepsilon_2,\ldots, \varepsilon_{k-1})}_{k-1})^2(1-\varepsilon_k)>1,
\end{equation}
and $\beta^{(\varepsilon_2,\ldots, \varepsilon_k)}_k>0$ satisfies
\begin{equation}
\label{eq:2.21N}
\beta^{(\varepsilon_2,\ldots, \varepsilon_k)}_k=1-e^{-q^{(\varepsilon_2,\ldots, \varepsilon_k)}_k\beta_k^{(\varepsilon_2,\ldots, \varepsilon_k)}}.
\end{equation}

Clearly, because of the deletions we have made and the use of
smaller connection probabilities, we have
\begin{equation}
\label{eq:2.22N}
Q_{\rm perc}\geq \beta_1 \prod^\infty_{k=1}\beta^{(\varepsilon_2,\ldots, \varepsilon_k)}_k
\end{equation}
(connections at different distances are independent).

Letting $\varepsilon_2, \varepsilon_k,\ldots \searrow 0$ (in this order), we have from (\ref{eq:2.20N}) and (\ref{eq:2.21N}): $q^{(\varepsilon_2)}_2 \nearrow c_2\beta^2_1$, hence $\beta^{(\varepsilon_2)}_2 \nearrow \beta_2$, where $\beta_2$ satisfies (\ref{eq:2.2}) for $k=2$, hence $q^{(\varepsilon_2, \varepsilon_3)}_3\nearrow c_3 \beta^2_2$, hence $\beta^{(\varepsilon_2, \varepsilon_3)}_3\nearrow \beta_3$, where $\beta_3$ satisfies (\ref{eq:2.2}) for $k=3$, and so on. Therefore,  we obtain (\ref{eq:2.7N}) from (\ref{eq:2.22N}).

Note that the cascade of giant components containing ${\mathbf 0}$ will always remain, and in any case it would not be deleted according to (\ref{eq:2.13N}), (\ref{eq:2.17N}), etc., a.a.s..

We now prove that
\begin{equation}
\label{eq:2.25'}
Q_{\rm perc}\leq \prod^{\infty}_{k=1}\beta_k.
\end{equation}
The approach is analogous to the one used for proving (\ref{eq:2.7N}), but  now we replace the random connection probabilities between giant components by upper bounds. Since non-giant components in $k$-balls have no effect on percolation, we consider only the connections between giant components.

For $k=2$, from (\ref{eq:2.9N}), (\ref{eq:2.11N}),
(\ref{eq:2.12N}), and (A) for the {\it ER} graph ${\mathcal G}(N,
c_1/N)$ we have
\begin{eqnarray*}
p^{N,2}_{ij}&\leq & \frac{c_2\beta^2_1}{N}+c_2
\biggl|\frac{|G^{(1)}_i||G^{(1)}_j|}{N^3}-\frac{\beta^2_1}{N}\biggr|\\
&=&\frac{c_2\beta^2_1}{N}+\frac{c_2}{N}\biggl|\frac{|G^{(1)}_i||G^{(1)}_j|}{N^2}-\beta^2_1\biggr|\\
&<& \frac{c_2}{N}(\beta^2_1+\delta_2)\quad {\rm a.a.s.}
\end{eqnarray*}
for any $\delta_2>0$. Now we connect the $1$-balls in a $2$-ball with  probability $c_2(\beta^2_1+\delta_2)/N$, and we consider the resulting level $2$ giant components $G^{(2)}_i$ in $2$-balls (we keep the notation $G^{(2)}_i$ used before, but now the connection probabilities are different).

For $k=3$, two giant components $G^{(2)}_i, G^{(2)}_j$, $i\neq j$, in a $3$-ball  are connected, by (\ref{eq:1.4}), as above, with probability
\begin{equation}
\label{eq:2.26'}
p^{N,3}_{ij}=1-\biggl(1-\frac{c_3}{N^5}\biggr)^{|G^{(2)}_i||G^{(2)}_j|},
\end{equation}
where $|G^{(2)}_i|=||G^{(2)}_i||N$, and, by (A),
\begin{equation}
\label{eq:2.27'}
||G^{(2)}_i||\sim \beta^{(\delta_2)}_2N\quad {\rm a.a.s.},
\end{equation}
where $\beta^{(\delta_2)}_2>0$ satisfies
\begin{equation}
\label{eq:2.28'}
\beta^{(\delta_2)}_2=1-e^{-q^{(\delta_2)}_2 \beta^{(\delta_2)}_2},
\end{equation}
with
\begin{equation}
\label{eq:2.29'}
q^{(\delta_2)}_2=c_2(\beta^2_1+\delta_2)
\end{equation}
Proceeding the same way as above, it follows from (\ref{eq:2.26'})-(\ref{eq:2.29'}) that
\begin{eqnarray*}
p^{N,3}_{ij}&\leq & \frac{c_2(\beta_2^{(\delta_2)})^2}{N} +c_3 \biggl|
\frac{|G^{(2)}_i||G^{(2)}_j|}{N^5}-\frac{(\beta_2^{(\delta_2)})^2}{N}\biggr|\\
&=&\frac{c_3(\beta_2^{(\delta_2)})^2}{N}+\frac{c_3}{N}\biggl|
\frac{||G^{(2)}_i||||G^{(2)}_j||}{N^2}-(\beta_2^{(\delta_2)})^2\biggr|\\
&<& \frac{c_3}{N}((\beta_2^{(\delta_2)})^2+\delta_3)\quad {\rm a.a.s.}
\end{eqnarray*}
for any $\delta_3>0$.

Iteration of this scheme leads to
\begin{equation}
\label{eq:2.30'}
Q_{\rm perc}\leq \beta_1 \prod^\infty_{k=1}\beta^{(\delta_2,\ldots, \delta_k)}_k,
\end{equation}
and $\beta^{(\delta_2,\ldots, \delta_k)}_k\searrow \beta_k$ as $\delta_2,\delta_3,\ldots \searrow 0$,
where $\beta_k$ satisfies (\ref{eq:2.2}).

We then obtain (\ref{eq:2.25'}) from (\ref{eq:2.30'}).

Since $Q_{{\rm perc}}$ refers to a special form of percolation, we
have   $Q_{\rm{perc}}\leq P_{\rm{perc}}$. In order to see that cascade
percolation is
equivalent to asymptotic percolation, and therefore $Q_{\rm{perc}}=P_{\rm{perc}}$,
which proves (\ref{eq:2.5}) by (\ref{eq:2.7'}), it remains only to show that, as we assumed at the beginning of the proof,
for each $k$ all $1$-step  connections from a $k$-ball are in the $(k+1)$-ball it belongs to a.a.s.. This  is done  later on in Lemma 2.5.

\hfill$\Box$

In the proof of Theorem 2.2 the random graph ${\mathcal G}_{N,k}$
has $N$ vertices which are the $(k-1)$-balls in a $k$-ball, and
each $(k-1)$-ball has an internal structure which contains the
giant and the other components in it.  We replaced the connection
probabilities between $(k-1)$-balls by the random connection
probabilities between their giant components because they
determine percolation, and then we replaced these random
connection probabilities by upper and lower (deterministic)
connection probabilities in order to obtain upper and lower bounds
for $Q_{\rm perc}$. The non-giant components have no part in
percolation, but still each $(k-1)$-ball contained all its
$N\;(k-2)$-balls (with their internal structures), and this was
relevant in the proof of the theorem (see Remark 2.4). The random
graphs ${\mathcal G}_{N,k}$ can be redefined so that their
vertices contain only giant components at previous levels, with
random connection probabilities between the vertices. In this
sense the graphs ${\mathcal G}_{N,k}$ are ``doubly stochastic''.
This way the {\it cascade of giant components in a $k$-ball} is
more clearly depicted: all the points of $\Omega_N$ in a vertex of
${\mathcal G}_{N,k}$ are connected, and for each $j=2,\ldots, k$,
the giant component at level $j$ consists of connected giant
components at level $j-1$. We denote by $C^{N,k}$ the size of the
cascade of giant components in a $k$-ball (i.e., the number of
points of $\Omega_N$ in the cascade). The growth of $C^{N,k}$ is
given in the following corollary. \vglue .25cm \noindent {\bf
Corollary 2.3} {\it
\begin{equation}
\label{eq:2.31'}
C^{N,k} \sim \biggl(\prod^k_{j=1}\beta_j\biggr)N^k\quad a.a.s.,
\end{equation}
where the $\beta_j$ satisfy (\ref{eq:2.2}).}
\vglue .15cm
\noindent
{\bf Proof.} For $k=1$ and a giant component $G^{(1)}$ of a $1$-ball, we have
\begin{equation}
\label{eq:2.32'}
|G^{(1)}|\sim \beta_1N\quad {\rm a.a.s.},
\end{equation}
by (\ref{eq:2.2}) and (A).

For $k=2$ and a level 2 giant component $G^{(2)}$ (with the true connection probabilities, not upper and lower bounds as before) in a $2$-ball, we have similarly, by (A), on the one hand,
\begin{equation}
\label{eq:2.33'}
||G^{(2)}||\gtrsim \beta^{(\varepsilon_2)}_2 N\quad {\rm a.a.s.},
\end{equation}
where $\beta^{(\varepsilon_2)}_2$ satisfies (\ref{eq:2.19'}), and on the other hand,
\begin{equation}
\label{eq:2.34'}
||G^{(2)}||\lesssim \beta^{(\delta_2)}_2 N\quad {\rm a.a.s.},
\end{equation}
where $\beta^{(\delta_2)}_2$ satisfies (\ref{eq:2.28'}). Since $\varepsilon_2$ and $\delta_2$ are both arbitrarily small, we have from (\ref{eq:2.33'}) and (\ref{eq:2.34'}),
\begin{equation}
\label{eq:2.35'}
||G^{(2)}||\sim \beta_2N\quad {\rm a.a.s..}
\end{equation}

Now, each vertex of $G^{(2)}$ contains $\sim \beta_1 N$ points a.a.s. (not all $N$ points, as in the proof of Theorem 2.2), therefore, from (\ref{eq:2.35'}),
$$
|G^{(2)}|\sim \beta_1\beta_2N^2\quad {\rm a.a.s..}
$$

By induction we obtain (\ref{eq:2.31'}). \hfill $\Box$ \vglue .5cm
\noindent {\bf Remark 2.4} We could consider the random graphs
${\mathcal G}_{N,k}$ with giant components alone in the proof of
Theorem 2.2, but then, instead of (\ref{eq:1.4}), which is
$c_k/N(N^{k-1})^2$, by Corollary 2.3 we would use connection
probabilities
$$
\frac{c_k}{N((\prod^{k-1}_{j=1}\beta_j)N^{k-1})^2}
$$
between points ${\mathbf x}$ and ${\mathbf y}$ such that $d({\mathbf x}, {\mathbf y})=k$, in order to
restrict the internal structures to the cascades of giant component in the $(k-1)$-balls containing
the points. The approach we have taken is more natural.

To complete the proof of Theorem 2.2 it remains to prove the following lemma.
\vglue.5cm
\noindent
{\bf Lemma 2.5} {\it For each $k\geq 0$, all external 1-step connections from a given $k$-ball are at distance $k+1$ from it a.a.s.}

\medskip

The proof of Lemma 2.5 is a special case of  the following
calculations, which are intended to  give some idea of what percolation paths might look
like as $N\to \infty$.

Let $Y^{(k)}$ denote the number of $1$-step connections from ${\mathbf 0}\in \Omega_N$ to points at
hierarchical distance $k\geq 1$ (or from any given point of $\Omega_N$ to points at distance $k$ from it;
 recall that $d(\cdot, \cdot)$ is translation-invariant). Since $Y^{(k)}$ is distributed
  Bin$(N^k-N^{k-1},c_k/N^{2k-1})$, then

\begin{equation}
\label{eq:2.36}
Y^{(1)}\Rightarrow {\rm Pois}(c_1)\quad{\rm as}\quad N\rightarrow \infty,
\end{equation}
and
\begin{equation}
\label{eq:2.22n}
Y^{(k)}\Rightarrow 0\quad {\rm as}\quad N\rightarrow \infty, \;\;k>1.
\end{equation}

In the case $k>1$,
\begin{equation}
\label{eq:2.23N}
P[Y^{(k)}>0]\sim \frac{c_k}{ N^{k-1}}
\quad {\rm as}\quad N\rightarrow \infty,
\end{equation}
and
\begin{equation}
\label{eq:2.24N}
P[Y^{(k)}=n|Y^{(k)}>0]\rightarrow \left\{\begin{array}{llc}
1, & n=1 \\
0, & n>1
\end{array}\right.{\rm as}\quad N\to \infty.
\end{equation}
Indeed, using $1-x \sim -\log x$ as $ x\rightarrow 1$,
\begin{eqnarray*}
P[Y^{(k)}>0] &=& 1-\biggl(1-\frac{c_k}{N^{2k-1}}\biggr)^{N^k-N^{k-1}}\\
&\sim & -(N^k -N^{k-1})\log \biggl(1-\frac{c_k}{N^{2k-1}}\biggr),
\end{eqnarray*}
and now using $x$log$(1-a/x)\sim -a$ as $x\rightarrow \infty$, $a>0$,
\begin{eqnarray*}
N^{k-1}P[Y^{(k)}>0] &\sim& -(N^{2k-1}-N^{2k-2})\log \biggl(1-\frac{c_k}{N^{2k-1}}\biggr)\\
&\sim & c_k,
\end{eqnarray*}
which proves (\ref{eq:2.23N}).
For the proof of (\ref{eq:2.24N}) we have
$$
P[Y^{(k)}=n|Y^{(k)}>0] = \frac{P[Y^{(k)}=n]}{P[Y^{(k)}>0]},
$$
and
\begin{eqnarray*}
P[Y^{(k)}=n]&=& \biggl({N^k-N^{k-1}\atop n}\biggr) \biggl(\frac{c_k}{N^{2k-1}}\biggr)^n \biggl(1-\frac{c_k}{N^{2k-1}}\biggr)^{N^k-N^{k-1}-n}\\
&\sim & \frac{(N^k-N^{k-1})!}{n!(N^k-N^{k-1}-n)!}\biggl(\frac{c_k}{N^{2k-1}}\biggr)^n,
\end{eqnarray*}
hence, by (\ref{eq:2.23N}),

\begin{eqnarray*}
\frac{P[Y^{(k)}=n]}{P[Y^{(k)}>0]}&\sim & \frac{(N^k-N^{k-1})!}{n!(N^k-N^{k-1}-n)!}\frac{N^{k-1}}{c_k} \biggl(\frac{c_k}{N^{2k-1}}\biggr)^n\\
&\rightarrow & \left\{
\begin{array}{cc}
1, & n=1\\
0, & n>1
\end{array}\right. {\rm as}\quad N\rightarrow \infty.
\end{eqnarray*}

The results (\ref{eq:2.36})-(\ref{eq:2.24N}) show that a.a.s.
almost all 1-step connections from ${\mathbf 0}$ will be at
distance 1 and the degree of ${\mathbf 0}$, and of any point, will be distributed
${\rm Pois}(c_1)$. Also, for each $k>1$ the probability of 1-step
connections from ${\mathbf 0}$  to points at distance $k$   is
negligible w.r.t. the probability of connections to points at
distance $k-1$, and given that it has connections at distance $k$,
there is only one such connection, as $N\to\infty$. In a similar
way it can be shown that the number of points in the giant
component of any  $1$-ball which have
external $1$-step connections at distance $2$ from it is asymptotically
Pois$(c_2\beta_1)$ as $N\to\infty$.

Similarly, let $Y^{(k,j)}$ denote the number of 1-step connections from a given $k$-ball
to $k$-balls at distance $j\geq k+1$ from it. $Y^{(k,j)}$ is distributed
Bin$(N^{j-k}-N^{j-1-k},p^{N,k,j})$, where
\begin{equation}
\label{eq:2.25N}
p^{N,k,j}=1-\displaystyle\left(1-\frac{c_j}{N^{2j-1}}\displaystyle\right)^{N^{2k}}
\end{equation}
is the probability of 1-step connections between two $k$-balls at distance $j$ from each other.

In the same way as above we have
\begin{equation}
\label{eq:2.26N}
p^{N,k,j}\sim \frac{c_j}{N^{2(j-k)-1}}\quad{\rm as}\quad N\to\infty,
\end{equation}
and  it follows that
\begin{eqnarray}
\label{eq:2.27N}
Y^{(k,k+1)}&\Rightarrow&{\rm Pois}(c_{k+1})\quad{\rm as}\quad N\to\infty,\\
\label{eq:2.28}
Y^{(k,j)}&\Rightarrow&0\quad{\rm as}\quad N\to\infty,\quad j>k+1,
\end{eqnarray}
and for $j>k+1$,
\begin{equation}
\label{eq:2.29}
P[Y^{(k,j)}>0]\sim \frac{c_j}{N^{j-k-1}}\quad{\rm as}\quad N\to \infty,
\end{equation}
and
\begin{equation}
\label{eq:2.30}
P[Y^{(k,j)}=n|Y^{(k,j)}>0]\to\left\{\begin{array}{ll}
1,&n=1\\
0,&n>1\end{array}\right. {\rm as}\quad N\to \infty.
\end{equation}

Analogous results can be obtained, using $(A)$, for the number of 1-step connections between a
 giant component in a $k$-ball and other such giant components at distance $j\geq k+1$ from it, with $c_j$ replaced by $c_j\beta^2_k$
 in (\ref{eq:2.26N}), (\ref{eq:2.27N}), (\ref{eq:2.29}). It can also be shown that the
 number of points in the cascade of giant components in  a $k$-ball with external $1$-step
 connections at distance $k+1$ is asymptotically
 Pois$(c_{k+1}\prod_{j=1}^k\beta_j)$ as $N\to\infty$.

From (\ref{eq:2.36})-(\ref{eq:2.24N}) and (\ref{eq:2.27N})-(\ref{eq:2.30}), we see that the same
pattern described above for connections between individual points
is repeated for connections between giant components of giant components  at the previous level in $k$-balls
for all $k$. In particular, for each $k$ all the external 1-step
connections  from the giant component in a given $k$-ball are at distance $k+1$ from it
a.a.s.

Note that these arguments also prove Lemma 2.5.

With these results one can also see that connections between $(k-1)$-balls in a $k$-ball through points outside the $k$-ball have probability $o(1/N)$, as stated in the Introduction.

These results suggest the following picture for the most likely form of cascade percolation paths:
 ${\mathbf 0}$ is in the giant component  of the 1-ball it belongs to. A Pois$(c_2\beta_1)$ number of points
 in this giant component have external single 1-step connections to points at distance 2 in such a way that the giant component is
 connected to the giant components of the other 1-balls, forming the level 2 giant component of the 2-ball
  they belong to. This  level $2$ giant component has a Pois $(c_3\beta_1\beta_2)$ number of points with external single $1$-step connections to other  level 2 giant components in the 3-ball they belong to, forming a  level 3 giant component. And so on along the cascade of giant components  for every distance $k$.  Each path starting from
${\mathbf 0}$ and going to infinity  this way is a percolation path.
\vglue .25cm
\noindent {\bf Example.} Let $c_k=a\log k$ for large $k$,
$a>0$. Then there is asymptotic percolation
if and only if $a>1$.

\vglue.5cm
\noindent
{\bf Remarks 2.6}

\noindent
{\bf 1.} Each of the following is a  sufficient condition for percolation:
$$
\sum\limits^\infty_{k=1}c^{-\delta}_k<\infty\quad\hbox{\rm  for some}\quad \delta>0,\;\; \liminf\limits_{k\rightarrow \infty} c_k/k>0,\;\; \liminf\limits_{k\rightarrow \infty} (c_{k+1}-c_k)>0.
$$

\noindent {\bf 2.} The degree sequence of $ER$ random graphs has
an approximate Poisson distribution.  On the other hand, the
random graphs of \cite{AB} have scale-free degree distributions. A
rigorous study of scale-free random graphs, including a coupling
with $ER$ random graphs, is contained in \cite{BR1, BR2}. We will
see that in the graphs ${\mathcal G}(N,p^{N,k})$  the degrees in a
neighborhood of $c_k$ have  an approximate scale-free behavior
with exponent $1/2$, for large $k$.

Recall that in the random graph ${\mathcal G}(N, p^{N, k})$ the
vertices are the $N\; (k-1)$-balls in a $k$-ball with $p^{N,k}$
given by (\ref{eq:1.5}). Neglecting the term $o(1/N)$ in
(\ref{eq:1.7}), we consider ${\mathcal G}(N,c_k/N)$.

Let $X_j$ denote the number of vertices in ${\mathcal G}(N,
c_k/N)$ of degree $j=0,1,\ldots$. Then the proportion $X_j/N$ has
an approximate Poisson distribution \cite{BR1} (Theorem 3), i.e.,
for each $j$ and any $0<\varepsilon <1$,
\begin{equation}
\label{eq:2.31}
P\biggl[(1-\varepsilon) \frac{c^j_k e^{-c_k}}{j!} \leq \frac{X_j}{N} \leq (1+\varepsilon )\frac{c^j_k e^{-c_k}}{j!}\biggr]\rightarrow 1 \;\;{\rm as}\;\;N\rightarrow \infty.
\end{equation}

Let $B_k =\{j:|j -c_k |<M\}$ for  some constant $M$. Recall that $c_k \rightarrow \infty$ as $k\rightarrow \infty$. Using the Stirling formula (in the form in \cite{F}, p.52), and the standard inequalities
$$
\biggl(1+\frac{x}{y}\biggr)^y <e^x <\biggl(1-\frac{x}{y}\biggr)^{-y},\quad 0<x<y, \;\;y>1,
$$
it is easy to show that for any $0<\delta<1$,
\begin{equation}
\label{eq:2.32}
(1-\delta )\frac{1}{\sqrt{2\pi}}j^{-1/2}\leq \frac{c^j_k e^{-c_k}}{j!} \leq \frac{1}{\sqrt{2\pi}}j^{-1/2}
\end{equation}
for all $j\in B_k$ and sufficiently large $k$.

Combining (\ref{eq:2.31}) and (\ref{eq:2.32}) we have, for any $0<\varepsilon <1$,
\begin{equation}
\label{eq:2.33}
P\biggl[(1-\varepsilon)\frac{j^{-1/2}}{\sqrt{2\pi}}\leq \frac{X_j}{N}\leq
(1+\varepsilon)\frac{j^{-1/2}}{\sqrt{2\pi}}\biggr]\rightarrow 1\;\;{\rm as}\;\;N\rightarrow \infty
\end{equation}
for  all $j\in B_k$ and sufficiently large $k$.

\section{Some open problems}
\label{sec:3}
\setcounter{equation}{0}
\subsection{A class of doubly stochastic random graphs}
\label{sub:3.1}
The essence of the proof of Theorem 2.2 is  that the random dependent connection probabilities $p^{N,k}_{ij}$ between giant components of $(k-1)$-balls in a $k$-ball   obey
\begin{equation}
\label{eq:3.1}
p^{N,k}_{ij}\sim \frac{c_k\beta^2_{k-1}}{N}\quad {\rm a.a.s.},\quad c_k\beta^2_{k-1}>1,
\end{equation}
in the sense that
\begin{equation}
\label{eq:3.2}
\frac{c_k\beta^2_{k-1}}{N}(1-\varepsilon)<p^{N,k}_{ij}<\frac{c_k \beta^2_{k-1}}{N}(1+\varepsilon)
\end{equation}
a.a.s. for each $k$ and arbitrary  $\varepsilon$ (such that
$c_k\beta^2_{k-1}(1-\varepsilon)>1$).  Therefore asymptotic
percolation in ${\mathcal G}_N$ is as if it took place along the
giant components of the {\it ER} random graphs ${\mathcal G}(N,
c_k\beta^2_{k-1}/N)$ instead of ${\mathcal G}(N, p^{N,k})$,
considering only  connections between $(k-1)$-balls.

This suggests the study of ``doubly stochastic'' random graphs of
the form ${\mathcal G}(N, \{p^N_{ij}\}_{i,j=1,\ldots, N})$ with
vertices $\{v_1,\ldots, v_N\}$ and connection probabilities
$p^N_{ij}$ for the vertices $v_i$ and $v_j$, where the $p^N_{ij}$
are random variables which may be dependent, such that
$p^N_{ij}\sim c/N, c>1$, for  all $ij$ a.a.s.. A natural question
is which results of the theory of {\it ER} random graphs
${\mathcal G}(N, c/N)$ can be extended to this more general class
of random graphs.

The usefulness of such extended results is illustrated in the following two problems, which could be treated if they were available.

\subsection{ Average distance in the cascade of giant components in a $k$-ball}
\label{sub:3.2}

Let $D(x,y)$ denote the (graph) distance between two points $x$ and $y$ (i.e. the length of the shortest path connecting the two points) chosen at
random in the cascade of giant components in a given $k$-ball, $k>1$. If the two points are at hierarchical distance $1$ from each other (i.e., they belong to the giant component of the same $1$-ball), then $D(x,y)$ grows like $\log N/\log c_1$ as  $N\rightarrow \infty$, according to Theorem 2.4.1 in \cite{D} and (1.1). If the hierarchical distance between the two points is larger, say $k$ (i.e., they are in different $(k-1)$-balls; note that this is the most likely outcome of the random choice of points, with probability $1-1/N$), then the expected result would be that $D(x,y)$ grows like
$$
\frac{(\log N)^k}{\prod^k_{j=1}\log (c_j \beta_{j-1}^2)}\quad {\rm as}\;N\rightarrow \infty,
$$
where the $\beta_{j-1}$ satisfy (\ref{eq:2.2}) and (\ref{eq:2.3}),
provided that  Theorem 2.4.1 in \cite{D} could be extended to the
random graphs with random connection probabilities $p^{N,k}_{ij}$
obeying (\ref{eq:3.1}).

Heuristically, the reasoning is as follows. Starting from a given point in  the giant component of a $1$-ball, the shortest path to one of the points  in the $1$-ball that have external connections at distance $2$ grows like  $\log N/\log c_1$ (these points are chosen at random in the giant component of the $1$-ball and all the shortest paths grow the same way). This happens in every $1$-ball. Next, in the giant component in the $2$-ball the
shortest path  from a given  level $1$ giant component to one of the  level $1$ giant components  that have connections at distance $3$ should grow like $\log N/\log (c_2\beta^2_1)$, according to (\ref{eq:3.1}), if the extended result holds. Hence the (graph) distance between two points chosen at random in the level $2$ giant component of  level $1$ giant components in a $2$-ball grows like the product $(\log N/c_1)(\log N/\log (c_2\beta^2_1))=(\log N)^2/\prod^2_{j=1}\log(c_j \beta^2_{j-1})$. Etc.

In contrast, the average distance in the giant component of the
{\it ER} random graph ${\mathcal G}(N^k, c/N^k), c>1$, grows like
$k\log N/\log c$. Hence, although the sizes of the giant clusters
in a $k$-ball of $\Omega_N$ (which has $N^k$ points)  and in the
{\it ER} graph ${\mathcal G}(N^k, c/N^k)$ are both of order $N^k$
(by (A) and Corollary 2.3), the typical length of paths is longer
in the hierarchical case (of order $(\log N)^k$). This happens
because two points chosen at random in the cascade in a $k$-ball
are most likely hierarchical distance $k$ apart, and the paths
connecting them have to go consecutively through $k$ hierarchical
levels (by cascade percolation).

\subsection{Central limit theorem for the size of the cascade of giant components in a $k$-ball}
\label{sub:3.3}

The following central limit theorem is proved in \cite{BBF} for
the fluctuation of the size $|G^N|$ of the giant component $G^N$
of a $ER$ random graph ${\mathcal G}(N,c/N), c>1:$
$$
\frac{|G^N|-\beta N}{N^{1/2}}\Rightarrow {\mathcal N}(0,
\sigma^2)\quad {\rm as}\quad N\rightarrow \infty,
$$
where $\Rightarrow$ denotes convergence in distribution, and
${\mathcal N}(0, \sigma^2)$ is the normal distribution with mean
$0$ and variance $\sigma^2=\beta (1-\beta)/\mu^2$,  where
$\beta>0$ satisfies $1-\beta=e^{-c\beta}$, and $-\mu$ is the slope
of $1-t-e^{-ct}$ at $t=\beta$ (see also \cite{D}, Theorem 2.5.3).

If this result could be extended to random graphs with random connection probabilities obeying (\ref{eq:3.1}), we would have the following central limit result. Recall that $C^{N,k}$ denotes the size of the  cascade of giant components in a  $k$-ball (see Corollary 2.3). Let
$$
X^{N,k}=\frac{C^{N,k}-(\prod^k_{j=1}\beta_j)N^k}{N^{k-1/2}},\quad k=1,2,\ldots.$$
Then
$$
X^{N,k}\Rightarrow {\mathcal N}\biggl(0,
\biggl(\prod^{k-1}_{j=1}\beta_j\biggr)^2 \sigma^2_k\biggr)
\quad{\rm as}\quad  N\rightarrow \infty
$$
for each $k$, where $\beta_k>0$ satisfies (\ref{eq:2.2}), $\sigma^2_k=\beta_k (1-\beta_k)/\mu^2_k$, and $-\mu_k$ is the slope of $1-t-e^{-c_k \beta^2_{k-1}t}$ at $t=\beta_k$.

We argue by induction. For $k=1$ the result is obtained directly
by the c.l.t. of \cite{BBF}. Now  we assume the result is true for
$k-1$ and consider the case $k$. Note that
$C^{N,k}=\sum^{||G^{N,k}||}_{j=1}C_j^{N,k-1}$, where $||G^{N,k}||$
is the size of the giant component of ${\mathcal G}(N, c_k
\beta^2_{k-1}/N)$, hence $||G^{N,k}||$ grows like $\beta_k N$. We
write $X^{N,k}$ as
\begin{eqnarray}
\label{eq:3.6.1}
X^{N,k} &=&\frac{1}{N^{k-1/2}}
\biggl(\sum^{||G^{N,k}||}_{j=1}C^{N,k-1}_j -\biggl(\prod^k_{j=1}\beta_j\biggr)N^k\biggr)\nonumber\\
&=&\frac{1}{N}\sum^{||G^{N,k}||}_{j=1}\frac{1}{N^{(k-1)-1/2}}\biggl(C^{N,k-1}_j
-\biggl(\prod^{k-1}_{i=1}\beta_i\biggr)N^{k-1}\biggr)\nonumber
\\&&\quad+\biggl(\prod^{k-1}_{i=1}\beta_i\biggr)\frac{||G^{N,k}||-\beta_kN}{N^{1/2}},
\end{eqnarray}
where $C^{N,k-1}_j$, $j=1,\ldots, ||G^{N,k}||$, are the sizes of
the cascades in the $(k-1)$-balls forming the giant component of
${\mathcal G}(N, c_k \beta^2_{k-1}/N)$. The $C^{N, k-1}_j$ are
i.i.d, and by the induction assumption,
\begin{equation}
\label{eq:3.4} \frac{1}{N^{(k-1)-1/2}}\biggl(C^{N, k-1}_j
-\biggl(\prod^{k-1}_{i=1}\beta_i\biggr)N^{k-1}\biggr)\Rightarrow
{\mathcal N}\biggl(0,
\biggl(\prod^{k-1}_{i=1}\beta_i\biggr)^2\sigma^2_{k-1}\biggr)
\end{equation}
for each $j$, as $N\rightarrow \infty$. It is easy to show  from (\ref{eq:3.4}) that the first term on the r.h.s. of (\ref{eq:3.6.1}) converges to $0$ in probability as $N\rightarrow \infty$. The result then follows from the second term in the r.h.s. of (\ref{eq:3.6.1}), if the extended central limit theorem for the random graph with connection probabilities $p^N_{ij}$ obeying (\ref{eq:3.1}) holds.

\vglue .25cm
\noindent
{\large\bf Appendix}
\vglue .15cm
\noindent
{\bf Proof of Lemma 2.1} (\ref{eq:2.2}) with $k=1$ and  $c_1> 2\log 2$ imply $\beta_1>1/2$. Then $c_2\beta_1^2>2$ log$2$, which by (\ref{eq:2.2}) with $k=2$ implies $\beta_2>1/2$, and so on, so $\beta_k >1/2$ for all $k$.

From (\ref{eq:2.2}) we have
\begin{eqnarray*}
1-\beta_k &=&\exp \{-c_k \beta^2_{k-1} (1-e^{-c_k \beta^2_{k-1}\beta_k})\}\\
&=&e^{-c_k\beta^2_{k-1}}\exp\biggl\{\frac{1}{\beta_k} c_k \beta^2_{k-1}
\beta_k e^{-c_k \beta^2_{k-1}\beta_k}\biggr\}\\
&&(\hbox{\rm using}\;\;\beta_k >1/2\;\;\hbox{\rm and}\;\;xe^{-x}\leq e^{-1})\\
&\leq & Ce^{-c_k\beta^2_{k-1}}.
\end{eqnarray*}
where $C=e^{2e^{-1}}$. Then, by (\ref{eq:2.2}) and $c_k \geq c_{k-1}$,
\begin{eqnarray*}
1-\beta_k & \leq & C \exp \{-c_{k-1}(1-e^{-c_{k-1}\beta^2_{k-2}\beta_{k-1}})^2\}\\
&\leq & C \exp \{-c_{k-1}(1-2e^{-c_{k-1}\beta^2_{k-2}\beta_{k-1}})\}\\
&\leq&Ce^{-c_{k-1}}\exp \{2c_{k-1}e^{-c_{k-1}/8}\}\\
&\leq & C_1 e^{-c_{k-1}},
\end{eqnarray*}
where $C_1=Ce^{16e^{-1}}$.

Hence $\sum e^{-c_k}<\infty$ implies $\sum (1-\beta_k)<\infty$, and therefore $\prod \beta_k>0$.

The reverse inequality is clear, since all $\beta_k<1$. \hfill $\Box$

\par\bigskip\noindent
{\bf Acknowledgement.} LGG thanks the Laboratory for Research in
Statistics and Probability, Carleton University, Ottawa, and the
Swiss Federal Institute of Technology (ETH), Z\"urich, for their
hospitality.

\bibliographystyle{amsplain}

\end{document}